\documentclass[12pt,a4paper]{article}

\usepackage{url}
\usepackage{defns}
\newcommand{\cJ}{{\cal J}}
\newcommand{\cD}{{\cal D}}
\newcommand{\uj}{u_j}
\newcommand{\sfrac}[2]{{\textstyle\frac{#1}{#2}}}
\newcommand{\ibc}{\textsc{ibc}}
\newcommand{\pde}{\textsc{pde}}

\begin{document}

\title{Holistic projection of initial conditions onto a finite 
difference approximation}

\author{A. J. Roberts 
\thanks{\protect\url{mailto:aroberts@usq.edu.au}}\\
\textit{Dept Maths \& Comput., University of Southern 
Queensland,}\\ \textit{Toowoomba, Queensland~4352, Australia}.
}
\maketitle

\begin{abstract}
Modern dynamical systems theory has previously had little to say about 
finite difference and finite element approximations of partial 
differential equations (\pde{}s) \cite{Archilla98}.
However, recently I have shown one way that centre manifold theory may 
be used to create and support the spatial discretisation of \pde{}s 
such as Burgers' equation \cite{Roberts98a} and the 
Kuramoto-Sivashinsky equation \cite{MacKenzie00a}.
In this paper the geometric view of a centre manifold is used to 
provide correct initial conditions for numerical discretisations 
\cite{Roberts97b}.
The derived projection of initial conditions follows from the physical 
processes expressed in the \pde{}s and so is appropriately 
conservative.
This rational approach increases the accuracy of forecasts made with 
finite difference models.
\end{abstract}



\section{Introduction}

Consider the equations for some physical field $u(x,t)$ evolving in 
space-time that we wish to model numerically.
Imagine a given initial field $u_0(x)$ and a finite difference model 
written in terms of $\uj(t)=u(x_j,t)$ for equi-spaced grid points 
$x_j=jh$ say; for example, in \S\ref{Scmth} for Burgers' 
equation~(\ref{eq:burg}) we find
\begin{equation}
     \frac{d\uj}{dt}+\frac{a}{2h}\mu\delta\, \uj^2\approx
    \frac{1}{h^2}\delta^2\uj 
    +\frac{a^2}{16}(\delta^2\uj^3-\uj^2\delta^2\uj) \,,
    \label{eq:dmod}
\end{equation}
in terms of the central difference operator $\delta\,\uj 
=u_{j+1/2}-u_{j-1/2}$ and central mean operator $\mu\,\uj 
=(u_{j+1/2}+u_{j-1/2})/2$\,.
One might expect that the correct initial condition for this 
discretisation is simply to project the initial field $u_0(x)$ onto 
the finite dimensional space of the model by setting the initial 
discretisation values to the value of the initial field at the grid: 
$\uj(0)=u_0(x_j)$.
But if the initial field is localised away from any grid point then 
physically we know to distribute the initial field among nearby grid 
points.
I use dynamical arguments to show that the correct initial condition 
is, to leading order, the correctly conservative element average
\begin{equation}
    \uj(0)\approx\frac{1}{h}\int_{x_j-h/2}^{x_j+h/2} u_0(x) \,dx\,.
    \label{eq:leadzj}
\end{equation}
This formula, and higher order corrections that involve neighbouring 
elements, are derived systematically herein.
For a numerical model, this is the first time a dynamical rationale 
has been used to provide initial conditions.

Such projection of initial fields onto the discretisation is supported 
by centre manifold theory \cite[e.g.]{Carr81}: the Relevance Theorem 
asserts that each of the nearby solutions of the governing \pde{} 
exponentially quickly in time approach a solution of the numerical 
model; this holds even for finite grid spacing~$h$.
The algebraic techniques developed by Roberts~\cite{Roberts97b}, based 
upon analysing with the aid of computer algebra the adjoint of a 
linearisation of the \pde{}, determines the initial condition for the 
discretisation so that we ensure the finite difference model 
faithfully tracks the correct particular solution of the \pde.

\section{Burgers' equation is discretised with centre manifold theory}
\label{Scmth}

Consider the dynamics of Burgers' equation
\begin{equation}
    u_t+auu_x=u_{xx}
    \label{eq:burg}
\end{equation}
as a prototype advection-diffusion equation.
Roberts~\cite{Roberts98a} first constructed finite difference 
approximations to the spatial derivatives using centre manifold theory 
to ensure nonlinear, subgrid-scale processes were systematically 
modelled.
We summarise the approach in this section.

Divide the spatial domain~$I$ into a number, say $m$, of 
elements of equi-size~$h$.
We analyse the dynamics of the elements away from any physical 
boundary to derive a discretisation for the interior of the domain.
Artificially crafted internal boundary conditions (\ibc's) between 
the elements are introduced:
\begin{equation}
    \left[u_x-\half au^2\right]=0\,,\quad
    (1-\gamma)h\overline{(u_x-\half au^2)}=\gamma\left[u\right]\,,
    \label{eq:ibc}
\end{equation}
where $[\phantom{u}]$ denote the jump across each internal boundary, 
$\overline{\phantom{u}}$ denotes the average value from the two sides 
of the boundary, and distinct from earlier work~\cite{Roberts98a} 
these \ibc's are expressed in terms of the flux $q=-u_x+\half au^2$\,.
See that when $\gamma=0$ the right-hand side of the second \ibc{} 
disappears so that the two conditions then completely insulate an 
element from its neighbours.
Whereas when $\gamma=1$, the left-hand side disappears and the two 
conditions ensure sufficient continuity of the physical field to 
recover Burgers' dynamics throughout the domain.

The centre manifold and the evolution thereon is straightforwardly 
constructed using the computer algebra algorithm described in 
\cite{Roberts98a,Roberts96a}.
Here we find the subgrid field in the $j$th~element is
\begin{eqnarray}
    && u(x,t) = \uj +\half ah\xi\uj^2 +\quarter a^2h^2\xi^2\uj^3
    +\gamma\left[\xi\mu\delta\,\uj 
    +\half\xi^2\delta^2\uj\right] \nonumber \\
     &  & {} +ah\gamma\left[ -{\textstyle\frac18}\xi( 
     \uj\delta^2\uj +\delta^2\uj^2 +4\uj^2) +{\textstyle\frac18}\xi^2( 
     2\uj\mu\delta\,\uj-\mu\delta\,\uj^2) 
     +\third\xi^2\uj\delta^2\uj  \right]  \nonumber \\
     &  & {} +a^2h^2\gamma\left[ 
     {\textstyle\frac{1}{16}}\xi( \uj\mu\delta\,\uj^2 +\mu\delta\uj^3) 
     -{\textstyle\frac{3}{32}}\xi^2( 3\uj^2\delta^2\uj +2\uj\delta^2\uj 
     -\delta^2\uj^3 +8\uj^3) 
     \right. \nonumber \\ &  & \left.\quad{}
     +{\textstyle\frac{1}{6}}\xi^3( 2\uj^2\mu\delta\,\uj 
     -\uj\mu\delta\,\uj^2) +{\textstyle\frac{5}{24}}\xi^4\uj^2\delta^2\uj 
     \right]  +\Ord{\gamma^2,a^3}\,,
     \label{eq:cm}
\end{eqnarray}
where $\xi=(x-x_j)/h$ ranges over $[-1/2,1/2]$\,.
The evolution on this centre manifold, when evaluated at $\gamma=1$ to 
restore continuity, forms the finite difference model~(\ref{eq:dmod}) 
for Burgers' equation: see that the first three terms 
in~(\ref{eq:dmod}) form a standard discretisation of each term but now 
appearing automatically from the discretisation when mediated by the 
flux form~(\ref{eq:ibc}) of the \ibc's; whereas the last term gives 
$\Ord{a^2}$ corrections to account for interactions between the 
nonlinear advection and the diffusive dissipation.
Such nonlinear modifications of standard discretisations can be 
extremely effective~\cite{Roberts98a}.

\begin{figure}[tbp]
    \centering
    \includegraphics[width=\textwidth]{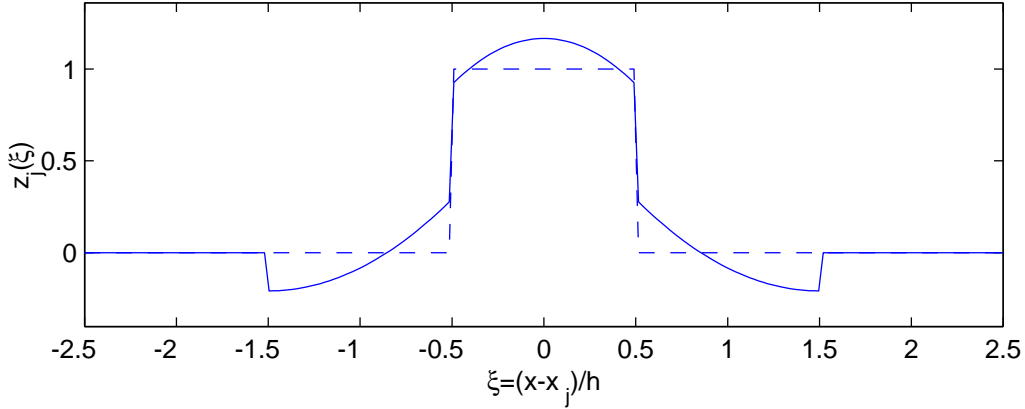}
    \caption{leading two orders of approximation to the projection 
    vectors~$z_j(x)$ for purely diffusive dynamics, errors: 
    $\Ord{\gamma}$, dashed; $\Ord{\gamma^2}$,  solid.  }
    \label{fig:zj}
\end{figure}

To find the correct initial condition, $\uj(0)$, for numerical models 
such as~(\ref{eq:dmod}) corresponding to any given field $u_0(x)$, we 
follow the procedure described in~\cite{Roberts97b}.
The aim is to determine projection vectors $z_j(x)$, such as those 
shown in Figure~\ref{fig:zj}, so that
\begin{equation}
    \left<z_j,u_0(x)-v(\vec u(0),x)\right>=0
    \quad\mbox{using}\quad
    \left<z,u\right>=\frac{1}{h}\int_I zu\,dx
    \label{eq:proj}
\end{equation}
as the inner product.
Now the dynamics linearised about the nonlinear centre manifold, 
$u=v(\vec u,x)$, is governed by the operator
\begin{displaymath}
    \cJ =\partial_x^2-av_x -av\partial_x\,,
\end{displaymath}
with  \ibc's linearised about~(\ref{eq:ibc}) of
\begin{equation}
    \left[u_x\right]=0\,,\quad
    (1-\gamma)h\overline{ (u_x- avu)}=\gamma\left[u\right]\,,
    \label{eq:libc}
\end{equation}
Then in the above inner product the adjoint of $\cJ$ is
\begin{equation}
    \cJ^{\dagger} z=\partial_x^2z+av\partial_xz\,,
    \quad\mbox{such that}\quad
         \left[z_x\right]=0\,,\quad
     (1-\gamma)h\overline{z_x}=\gamma\left[z\right]\,.
    \label{eq:adj}
\end{equation}
To find the projection vectors~$z_j(x)$ we start with the leading 
approximation $z_j(x)\approx\chi_j(x)$ corresponding 
to~(\ref{eq:leadzj}) and plotted in Figure~\ref{fig:zj}, where 
$\chi_j(x)$ denotes the characteristic function that is~$1$ in the 
$j$th~element and otherwise is~$0$.
Then successive corrections are sought by iteration to ultimately 
satisfy the appropriate version of the equations derived 
in~\cite{Roberts97b}: defining the dual operator $\cD z=\D 
tz+\cJ^{\dagger}z$ we must solve
\begin{equation}
    \cD z_j-\sum_i\left<\cD z_j,e_i\right>z_i=0\,,
    \label{eq:projdz}
\end{equation}
subject to the \ibc's in~(\ref{eq:adj}) and the normalisation 
condition
\begin{equation}
    \left<z_j,e_i\right>=\delta_{i,j}\,,
    \label{eq:norm}
\end{equation}
where $e_j=\partial v/\partial \uj$ is the tangent vector of the 
centre manifold.
We seek solutions in a power series in $\gamma$ to errors 
$\Ord{\gamma^{\ell}}$ corresponding to the finite difference 
approximation of stencil width $2\ell-1$.
A computer algebra program available from the author does all the 
necessary algebra.

\section{Project onto Burgers' discretisation}

In this section we solve to quantities with errors 
$\Ord{a^3,\gamma^2}$: the finite difference model for Burgers' 
equation is then~(\ref{eq:dmod}); and the corresponding centre 
manifold over the whole domain is given by~(\ref{eq:cm}).
Calculating to errors $\Ord{a^3,\gamma^2}$ the projection onto the 
numerical model must be orthogonal to
\begin{eqnarray}
    z_j&\!\approx\!&\left[1-\frac{h^2a^2}{16}\uj^2\right]\chi_j
    \nonumber\\&&{}
    +\gamma\left[ \left(\sfrac16-\xi^2\right)\chi_j
    +\left(-\sfrac1{12}+\half\xi+\half\xi^2\right)\chi_{j-1}
    +\left(-\sfrac1{12}-\half\xi+\half\xi^2\right)\chi_{j+1} \right]
    \nonumber\\&&{}
    +\frac{ha\gamma}{48}\left[ \left( -(12\xi-16\xi^3)\uj 
    +u_{j+1}-u_{j-1} \right)\chi_j
    \right.\nonumber\\&&\quad\left.{}
    +\left( +\uj-(3-6\xi+8\xi^3)u_{j-1} \right)\chi_{j-1}
    \right.\nonumber\\&&\quad\left.{}
    +\left( -\uj+(3+6\xi-8\xi^3)u_{j+1} \right)\chi_{j+1} \right]
    \nonumber\\&&{}
    +\frac{h^2a^2\gamma}{384}\left[ \left( 8(1+3\xi^2)\uj^2 
    +4\uj(u_{j+1}+u_{j-1}) +2(u_{j+1}^2+u_{j-1}^2) \right)\chi_j 
    \right.\nonumber\\&&\quad\left.{}
    +\left( 3\uj^2+4\uj u_{j-1}-(5+12\xi^2+16\xi^3)u_{j-1}^2 
    \right)\chi_{j-1}
    \right.\nonumber\\&&\quad\left.{}
    +\left( 3\uj^2+4\uj u_{j+1}-(5+12\xi^2-16\xi^3)u_{j+1}^2 
    \right)\chi_{j+1}
    \right]\,.
    \label{eq:zj1ad}
\end{eqnarray}
Higher order expressions may be straightforwardly computed by computer 
algebra.
I conclude by further interpreting the physical effects incorporated 
in the projection defined by the above~$z_j$.

\subsection{Linear diffusion}

\begin{figure}[tbp]
    \centering
    \includegraphics{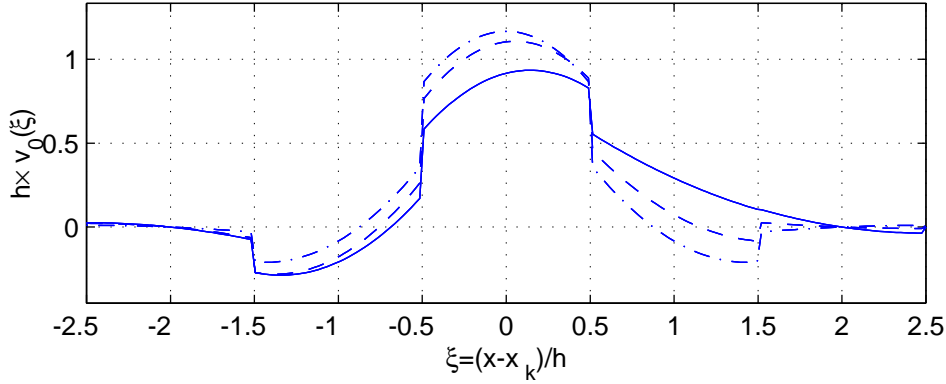}
    \caption{initial fields $u=v(\vec u(0),x)$ corresponding to a 
    unit-mass point release at: $\xi=0$, dot-dash; $\xi=1/4$, dashed; 
    $\xi=1/2$, solid.}
    \label{fig:ifld}
\end{figure}

Set $a=0$ in this subsection to analyse the linear diffusion equation 
$u_t=u_{xx}$\,.
Then the projection vector~(\ref{eq:zj1ad}), evaluated at $\gamma=1$ 
to recover the physically relevant case as plotted in 
Figure~\ref{fig:zj}, is
\begin{equation}
    z_j\approx\left(\sfrac76-\xi^2\right)\chi_j
    +\left(-\sfrac1{12}+\half\xi+\half\xi^2\right)\chi_{j-1}
    +\left(-\sfrac1{12}-\half\xi+\half\xi^2\right)\chi_{j+1}\,.
    \label{eq:zj1}
\end{equation}
To find the correct initial condition using this in~(\ref{eq:proj}) 
note that in these linear diffusion dynamics $\left<z_j,v(\vec 
u,x)\right>=\uj$ by the normalisation~(\ref{eq:norm}); thus here 
$\uj(0) =\left< z_j,u_0(x) \right>$\,.
For example, see that a point release in the $k$th element, 
$u_0(x)=\delta(x-x_k-h\eta)$, requires the slightly distributed 
initial condition
\begin{equation}
    h\uj(0)=\left(\sfrac{7}{6}-\eta^2\right)\delta_{k,j}
    +\left(-\sfrac{1}{12}-\half\eta+\half\eta^2\right)\delta_{k-1,j}
    +\left(-\sfrac{1}{12}+\half\eta+\half\eta^2\right)\delta_{k+1,j} 
    \,.
    \label{eq:diffptg}
\end{equation}
Such a specific initial condition corresponds via~(\ref{eq:cm}) to a 
field on the centre manifold as shown in Figure~\ref{fig:ifld} for the 
three cases
$\eta=0$, $1/4$ and~$1/2$.

See that these initial conditions ensure that the first moment of the 
numerical solution is correct for all time: in the numerical 
model~(\ref{eq:dmod}) the first moment is constant in time so it is 
enough to check that the first moment is correct in the initial 
conditions.
Define $\left<u\right>=\left<1,u\right>$\,, then for all time 
$\left<u\right>=1$ both in the model and in the exact solutions.
The first moment in the exact solution is its initial value 
$m_1=\left<(x-x_k)u_0(x)\right>=h\eta$\,; from~(\ref{eq:diffptg}) the 
first moment in the numerical model is the same
\begin{displaymath}
    m_1=\left<(x-x_k)v(\vec u(0),x)\right>=h\eta\,.
\end{displaymath}
However, the second moment $m_2=\left<(x-x_k)^2u\right>$ has 
$\Ord{h^2}$ errors: it evolves in time at the correct rate 
$dm_2/dt=2$, but the initial value is $h^2(\eta^2-1/6)$ instead of~$0$.
Determining the projection of initial conditions to higher orders in 
the coupling parameter~$\gamma$ obtains such higher order moments 
correctly.
Note that the rational approach adopted here does better than the 
usually chosen initial conditions which incur $\Ord{h}$ errors.

\subsection{Nonlinear dynamics}

Consider the $\Ord{a}$ terms from~(\ref{eq:zj1ad}) that 
modify~(\ref{eq:zj1}), namely
\begin{eqnarray*}
z_j'&=&{}
    +\frac{ha}{48}\left[ \left( (-12\xi+16\xi^3)\uj 
    +u_{j+1}-u_{j-1} \right)\chi_j
    \right.\nonumber\\&&\quad\left.{}
    +\left( +\uj+(-3+6\xi-8\xi^3)u_{j-1} \right)\chi_{j-1}
    \right.\nonumber\\&&\quad\left.{}
    +\left( -\uj+(+3+6\xi-8\xi^3)u_{j+1} \right)\chi_{j+1} \right]\,.
\end{eqnarray*}
Realise that the leading order effect of including these terms is to 
\emph{modify} the initial condition by $\left<z_j',u_0(x)\right>$\,.  
For example, if the initial field is approximately constant, 
$u_0(x)\approx U$, then 
\begin{displaymath}
    z_j'=\frac{Uha}{48}\left[ 
    (6\xi-8\xi^3)(\chi_{j+1}-2\chi_j+\chi_{j-1}) 
    +2(\chi_{j+1}-\chi_{j-1}) \right]\,;
\end{displaymath}
that the coefficients of the characteristic functions~$\chi_k$ sum to 
zero reflects that the the projection conserves the field~$u$.
More specifically, if $u_0(x)$ is~$U$ except for a symmetric bump in 
the $k$th element, then as well as the direct symmetric distribution 
identified for linear diffusion, the component in 
$\chi_{j+1}-\chi_{j-1}$ causes $u_{k-1}(0)$ to increase and 
$u_{k+1}(0)$ to decrease by an amount proportional to~$Ua$ reflecting 
that the self advection of the bump is not as great as that induced by 
assigning the mass of the bump solely to~$u_k(0)$.
Conversely, for an antisymmetric perturbation in the $k$th element, 
positive to the left of $x_k$, the component in 
$\chi_{j+1}-2\chi_j+\chi_{j-1}$ increases $u_k(0)$ and decreases 
$u_{k\pm 1}(0)$ in proportion to~$Ua$ to reflect the increased delay 
in~$u$ advecting out of the $k$th element because more of it is 
further to the left initially.
The $\Ord{a^2}$ terms in~(\ref{eq:zj1ad}) reflect more subtle physical 
processes.

\section{Conclusion}

Based upon the method of analysis and the discussion in the previous 
sections, we deduce that this centre manifold approach to finding 
correct initial conditions for finite difference models accounts for 
subgrid scale processes that occur as initial transients decay.
No other method does this.

Extensions of this approach to higher spatial dimensions is 
straightforward.
For example, consider the class of diffusive \pde's
\begin{displaymath}
	 \D tu=\delsq u+f(u,\grad u)\,,
\end{displaymath}
where $f$ represents nonlinear reaction or advection effects.
After tessellating space into finite elements---using \ibc's of the 
form (cf (\ref{eq:ibc}))
\begin{displaymath}
    \left[q_n\right]=0
    \quad\mbox{and}\quad
    (1-\gamma)h\overline{q_n}=\gamma[u]
\end{displaymath}
where $q_n$ is the flux of~$u$ normal to the internal boundary and $h$ 
is a size of the element---the fundamental problem in constructing a 
model is simply to solve Poisson's equation with forced Neumann 
boundary conditions on each element.
The adjoint of this problem lies at the heart of the 
dual~(\ref{eq:projdz}) for determining initial conditions of the 
approximation.
Although these sub-grid problem may itself need to be done 
numerically, in the simplest case of a regular tessellation it need 
only be done once for each term in the model, just like the 
computation of the interaction terms in a traditional finite element 
approximation.

In the case where there are variations in the size or shape of the 
elements of the discretisation, one would build formulae for the 
approximation parametrised by the shapes of the element and those 
neighbouring elements to which it is coupled by the \ibc's.
The algebraic detail becomes more complicated but the principles are 
the same.

\small
\bibliography{ajr,new,bib}
\bibliographystyle{unsrt}

\end{document}